\documentclass[11pt]{amsart}

\usepackage{amssymb}
\usepackage{euscript}
\usepackage{xspace}

\newcommand{\ie}{i.e.,\xspace}
\newcommand{\eg}{e.g.,\xspace}

\newtheorem{thm}{Theorem}
\newtheorem{prop}{Proposition}
\newtheorem{lemma}{Lemma}
\newtheorem{cor}{Corollary}
\newenvironment{conj}{\par\vspace{.5em}
        \noindent\textbf{Conjecture.}\em}{\em\par\vspace{.5em}}

\theoremstyle{definition}
\newtheorem{defn}{Definition}

\theoremstyle{remark}
\newtheorem{rem}{Remark}

\newcommand{\cal}{\EuScript}
\newcommand{\M}{\cal M}
\newcommand{\C}{\cal C}

\newcommand{\Sn}{\mathbf{S}}
\newcommand{\Tn}{\mathbf{T}}

\renewcommand{\leq}{\leqslant}
\renewcommand{\geq}{\geqslant}
\renewcommand{\kappa}{\varkappa}
\newcommand{\D}{\partial}

\newcommand{\R}{\mathbb R}
\renewcommand{\S}{\mathbb S}

\newcommand{\norm}[1]{\left|#1\right|}
\newcommand{\Norm}[1]{\|#1\|}

\newcommand{\tr}{\operatorname{tr}}

\newcommand{\diag}{\operatorname{diag}}

\newcommand{\dist}{\operatorname{dist}}

\newcommand{\compose}{\operatorname{\scriptstyle\circ}}

\newcommand{\N}{{\mathbf N}}
\newcommand{\X}{{\mathbf X}}
\newcommand{\x}{{\mathbf x}}
\renewcommand{\v}{{\mathbf v}}
\newcommand{\e}{{\mathbf e}}
\newcommand{\0}{{\mathbf 0}}

\newcommand{\alphab}{{\boldsymbol{\alpha}}}
\newcommand{\betab}{{\boldsymbol{\beta}}}
\newcommand{\gammab}{{\boldsymbol{\gamma}}}
\newcommand{\lambdab}{{\boldsymbol{\lambda}}}
\newcommand{\mub}{{\boldsymbol{\mu}}}
\newcommand{\Gammab}{{\boldsymbol{\Gamma}}}

\newcommand{\aut}[1]{{\sc #1}}
\newcommand{\tit}[1]{{\em #1\/}}
\newcommand{\vol}[1]{{\bf #1}}
\newcommand{\yr}[1]{(#1)}
\newcommand{\pp}[2]{#1--#2}

\newcommand{\st}{\colon\>}

\begin{document}

\title[Isometric Embeddings]{A priori bounds for \\
        co-dimension one isometric embeddings}
\author{Yanyan Li}
\address{Department of Mathematics\\ Rutgers University\\
        New Brunswick, New Jersey 08903}
\email{yyli@math.rutgers.edu}
\thanks{The research of the first author
        was supported in part by
        the Alfred P. Sloan Foundation Research Fellowship,
         NSF grant DMS-9401815 and DMS-9706887,
        and a Rutgers University Research Council grant.\\ \indent
        The research of the second author
        was supported in part by NSF Grant~DMS-9404523 and
        DMS~9704760.}
\author{Gilbert Weinstein}
\address{Department of Mathematics\\
        University of Alabama at Birmingham\\
        Birmingham, Alabama 35205}
\email{weinstei@math.uab.edu}

\date{\today}

\begin{abstract}
Let  $\X\colon(\S^n,g)\to\R^{n+1}$ be a $\C^4$ isometric embedding of a
$\C^4$ metric $g$ of non-negative sectional curvature on $\S^n$ into the
Euclidean space $\R^{n+1}$.  We prove a priori bounds for the trace of
the second fundamental form $H$, in terms of the scalar curvature $R$ of
$g$, and the diameter $d$ of the space $(\S^n,g)$.  These estimates give
a bound on the extrinsic geometry in terms of intrinsic quantities. They
generalize estimates originally obtained by Weyl for the case $n=2$ and
positive curvature, and then by P. Guan and the first author for
non-negative curvature and $n=2$.  Using $\C^{2,\alpha}$ interior
estimates of Evans and Krylov for concave fully nonlinear elliptic
partial differential equations, these bounds allow us to obtain the
following convergence theorem: For any $\epsilon>0$, the set of metrics
of non-negative sectional curvature and scalar curvature bounded below
by $\epsilon$ which are isometrically embedable in Euclidean space
$\R^{n+1}$ is closed in the H\"older space $\C^{4,\alpha}$,
$0<\alpha<1$. These results are obtained in an effort to understand the
following higher dimensional version  of the Weyl embedding problem
which we propose: \emph{Suppose that $g$ is a smooth metric of
non-negative sectional curvature and positive scalar curvature on $\S^n$
which is locally isometrically embeddable in $\R^{n+1}$.  Does
$(\S^n,g)$ then admit a smooth global isometric embedding
$\X\colon(\S^n,g)\to\R^{n+1}$?}
\end{abstract}
\maketitle

\section{Introduction}

In 1916, H.~Weyl posed the following problem: \emph{Given a metric $g$
of positive Gauss curvature on the sphere $\S^2$, is there an embedding
$\X\colon\S^2\to\R^3$ such that the metric induced on $\S^2$ by this
embedding is $g$?} Such an embedding $\X\colon (\S^2,g) \to \R^3$ is
called {\em isometric}, and satisfies the following system of nonlinear
partial differential equations:
\begin{equation}        \label{eq:iso}
        \nabla_i X \cdot \nabla _j X = g_{ij}.
\end{equation}

In~\cite{weyl16}, Weyl suggested the continuity method to attack the
problem and obtained a priori estimates up to the second derivatives of
the embedding.  The a priori estimate of the second derivative is a
consequence of the Weyl inequality which gives a bound of the mean
curvature by intrinsic quantities for strictly convex closed surfaces;
see Theorem~\ref{thm:weyl}. The main obstacle to the solution was the
lack of $\C^{2,\alpha}$ a priori estimates for the embedding.

Later H.~Lewy solved the problem under the assumption that the metric
$g$ is real analytic; see~\cite{lewy38}. It is interesting to point out
that Lewy did not use Weyl's a priori estimate. In~\cite{nirenberg53a},
L.~Nirenberg gave a beautiful proof for any metric $g$ of class $\C^4$.
He established, among other things, the $\C^{2,\alpha}$ a priori
estimate of the embedding using strong a priori estimates  he derived
earlier for solutions of  fully non-linear elliptic partial differential
equations in two variables; see~\cite{nirenberg53}. An entirely
different approach was taken independently by A.~D.~Alexandroff and
A.~V.~Pogorelov; see~\cite{alexandroff48,pogorelov49,pogorelov49a}. The
Weyl estimate was later generalized to the case of non-negative
curvature by P.~Guan and Y.~Li in~\cite{guanli94}.  From their estimate,
they obtained a $\C^{1,1}$ embedding result for metrics of non-negative
Gauss curvature; see also~\cite{hong95} for a different approach to the
$\C^{1,1}$ embedding result.

The main result of this paper is a Weyl-type estimate, see
Theorem~\ref{thm:guanli94} below, which generalizes to higher dimensions
the estimate of P.~Guan and the first author. The theorem asserts that
for any convex closed hypersurface in $\R^{n+1}$, ($n\ge 2$), one can
bound the mean curvature $H$ in terms of the scalar curvature of the
induced metric $g$, its Laplacian, and the diameter.

Denote by $\M^k(\S^n)$ the space of metrics $g$ of non-negative
sectional curvature on $\S^n$ which have $k$ continuous derivatives, and
by $\M^k_+(\S^n)\subset\M^k(\S^n)$ the subset consisting of those
metrics which have positive sectional curvature.  Similarly, if
$0<\alpha<1$, we denote by $\M^{k,\alpha}(\S^n)\subset\M^k(\S^n)$ the
space of metrics $g\in\M^k(\S^n)$ whose $k$-th derivatives are H\"older
continuous, and by $\M^{k,\alpha}_+(\S^n)\subset\M^{k,\alpha}(\S^n)$ the
subset consisting of those metrics which have positive sectional
curvature.  For convenience, we will denote $\M^{k,0}(\S^n)=\M^k(\S^n)$,
$\M^{k,0}_+(\S^n)=\M^k_+(\S^n)$, and adopt the same convention for
spaces of functions, \ie $\C^{k,0}=\C^k$.

\begin{thm}     \label{thm:guanli94}
Let $g\in\M^4(\S^n)$, and let $\X\colon(\S^n,g)\to\R^{n+1}$ be a $\C^4$
isometric embedding.  Let $H$ be the trace of the second fundamental
form of $\X$, and let $R$ be the scalar curvature of $g$.  Then the
following inequality holds
\begin{equation}        \label{eq:guanli}
        H^2 \leq C d^2 \,
        \sup_{\S^n}\left( 2R^2 - \Delta R + \frac{(n-1)^2}{64 d^2} R
        \right),
\end{equation}
where $C= 4(n-1)^{-2}e^{(n-1)/4}$, and $d$ is the diameter of
$(\S^n,g)$.
\end{thm}

\begin{rem}\ The above bound of $H$ involves four derivatives of $g$.
It is not known whether it is possible to bound $H$ by quantities
involving only up to third derivatives of $g$.  If the sectional
curvature of $g$ is strictly positive and $n\geq3$, then we can bound
$H$ by quantities involving only up to second derivatives of $g$. See
Remark~\ref{rem:est} below and Theorem~\ref{thm:c2bound} for details.
\end{rem}

As the first step in establishing Theorem~\ref{thm:guanli94}, we show
that for any convex closed surface with positive scalar curvature in
$\R^{n+1}$ ($n\ge 2$), one can bound the mean curvature $H$ in terms of
the scalar curvature of the induced metric $g$  and its Laplacian.  This
is a direct generalization of the Weyl estimate.

\begin{thm}     \label{thm:weyl}
Let $g\in\M^4(\S^n)$, and let $\X\colon(\S^n,g)\to\R^{n+1}$
be a $\C^4$ isometric embedding.  Let $H$
be the trace of the second fundamental form of $\X$, and let $R$ be the
scalar curvature of $g$.  Suppose that $R>0$.
Then the following inequality holds:
\begin{equation}        \label{eq:weyl}
        H^2 \leq \sup_{\S^n}\left( 2R - \frac{1}{R}\, \Delta R\right).
\end{equation}
\end{thm}

\begin{rem} \label{rem:est}
An estimate similar to~\eqref{eq:weyl} under the stronger hypothesis
that $g$ has positive sectional curvature was established
in~\cite{yau80}. In Theorem~\ref{thm:c2bound} we show that when $g$ has
positive sectional curvature and the dimension $n\ge 3$, $H$ can be
bounded in terms only of the lower bound of the sectional curvature and
the upper bound of the Ricci curvature.
\end{rem}

Using $\C^{2,\alpha}$ interior estimates of Evans and Krylov for concave
fully nonlinear elliptic partial differential equations,
see~\cite{evans82,krylov82,krylov84,caffarelli95}, the a priori bound in
Theorem~\ref{thm:weyl} allows us to obtain the following convergence
theorem.

\begin{thm}     \label{thm:convergence}
For any $\epsilon>0$, the set of metrics of non-negative sectional
curvature and scalar curvature bounded below by $\epsilon$ which are
isometrically embedable in Euclidean space $\R^{n+1}$ is closed in the
H\"older space $\C^{4,\alpha}$, $0<\alpha<1$.
\end{thm}

We note that $\C^{2,\alpha}$ estimates up to the boundary for concave
fully nonlinear elliptic partial differential equations were
independently established by Caffarelly, Nirenberg, and
Spruck~\cite{cns}, and Krylov~\cite{krylov84}.

In~\cite[Problem 53]{yau82}, S.~T.~Yau posed the following problem:
\emph{Can one generalize Weyl's embedding problem to higher dimensions?
More precisely, given a compact $n$-dimensional Riemannian manifold
$(M,g)$ of positive sectional curvature, is there an isometric immersion
$\X\colon (M,g)\to \R^{n(n+1)/2}$?}  The Cartan-Janet dimension
$n(n+1)/2$ is the smallest so that formally, \eg in the analytic class,
the problem has a local solution, regardless of curvature;
see~\cite{jacobowitz82}. In this direction, see
also~\cite{berger83,bryant83}, and the references therein.

Here, we wish to consider a different generalization of the Weyl
embedding problem: \emph{Can one give an intrinsic characterization of
those metrics $g$ with non-negative sectional curvatures on $\S^n$ for
which there is a co-dimension one isometric embedding
$\X\colon(\S^n,g)\to\R^{n+1}$?} Note that a necessary condition for such
an embedding to exist is that $g$ be locally isometrically embeddable in
$\R^{n+1}$.  In dimension $n\geq3$, this is a non-trivial restriction.
We formulate the following conjecture

\begin{conj}\  \label{conj:embed}
Let  $g$ be a smooth metric of non-negative sectional curvature and
positive scalar curvature on $\S^n$ which is locally isometrically
embeddable in $\R^{n+1}$.  Then  $(\S^n,g)$  admits a smooth global
isometric embedding $\X\colon(\S^n,g)\to\R^{n+1}$.
\end{conj}

The estimates we prove here are obtained in an effort to confirm the
conjecture. When $n\ge 3$ and $g$ has positive sectional curvature, any
local isometric embedding is known to be rigid, see
Section~\ref{sec:gauss}. Therefore the existence of a global immersion
follows by a standard monodromy argument from local embeddability and the fact
that $\S^n$ is simply connected. The immersion has to be
an embedding due to a theorem of Hardamard, see
~\cite[Theorem 2.11
on page 94]{spivak75} for a proof in the $2$-dimensional case.
The argument actually applies in any dimension. Thus, in view of
Theorem~\ref{thm:convergence}, the conjecture would be confirmed if one
could approximate any metric $g$ of non-negative curvature and positive
scalar curvature which is locally isometrically embeddable in
$\Bbb R^{n+1}$, by a sequence $g_i$ of metrics of positive curvature which
are also locally isometrically embeddable in $\R^{n+1}$. Conversely,
assuming the conjecture holds, any locally embeddable metric of
non-negative curvature and positive scalar curvature can be approximated
by embeddable metrics of positive curvature as follows from the
following simple proposition:

\begin{prop} \label{prop:approx}
Let $g\in\M^{k,\alpha}(\S^n)$ be isometrically embeddable in $\R^{n+1}$,
then there is a sequence $g_i\in\M^{k,\alpha}_+(\S^n)$ isometrically
embeddable in $\R^{n+1}$which converges to $g$.
\end{prop}

\begin{proof}
Let $\X\colon(\S^n,g)\to\R^{n+1}$ be an isometric embedding.  We can
write $\X$ as a graph over $\S^n$:
\[
        \X=\rho\x,\quad \x\in\S^n,\quad \rho=\rho(\x)>0.
\]
Then the first and second fundamental forms are given by:
\begin{align*}
        g_{ij} &= \rho^2 \gamma_{ij} + \rho_i\rho_j\\
        \chi_{ij} &= \rho^2\gamma_{ij} + 2\rho_i\rho_j - \rho\rho_{ij},
\end{align*}
where $\gamma_{ij}$ is the standard metric on $\S^n$, and the subscripts
denote covariant derivatives with respect to $\gamma$,
see~\cite{guanspruck93}.  If we substitute $\rho=u^{-1}$, then we have
for the second fundamental form:
\[
        \chi_{ij} = u^{-3} (u\gamma_{ij} + u_{ij}).
\]
Since $g$ has non-negative curvature, $\chi_{ij}$ is positive
semi-definite, see Section~\ref{sec:weyl}.  Hence, we have:
\[
        (u\gamma_{ij} + u_{ij}) \geq 0.
\]
If we set $u^\epsilon=u+\epsilon$, then we get:
\[
        \chi^\epsilon_{ij} = (u^\epsilon)^{-3}
        (u\gamma_{ij} + u_{ij} + \epsilon\gamma_{ij}) > 0,
\]
Hence the first fundamental forms $g^\epsilon_{ij}$ of
$\X^\epsilon=(u^\epsilon)^{-1}\x$ have positive curvature.  Clearly,
$g^{\epsilon}_{ij}$ converge to $g$ as $\epsilon\to0$.
\end{proof}

The paper is organized as follows. First, in Section~\ref{sec:weyl}, we
prove Theorem~\ref{thm:weyl}. In Section~\ref{sec:converge}, we prove
Theorem~\ref{thm:convergence}.  Then, in Section~\ref{sec:weyltype}, we
prove Theorem~\ref{thm:guanli94}. Finally, in Section~\ref{sec:gauss},
we show that the once-contracted Gauss equations can be solved for the
second fundamental form, when the sectional curvature of $g$ is
positive. As a corollary, we obtain a priori bounds on the second
fundamental form and hence also on the second derivatives of the
embedding, which depend only on two derivatives of $g$. This also
provides, when $n=3$, a local explicit criterion for any metric $g$ of
positive sectional curvature on $\S^3$ to be locally isometrically
embeddable in $\R^4$.

\section{Weyl Estimates in Higher Dimension}
\label{sec:weyl}

In this section, we will prove Theorem ~\ref{thm:weyl}, a higher
dimensional analogue of the Weyl estimate~\cite{weyl16}.  However, we
first note that Theorem~\ref{thm:weyl} allows us to establish a priori
bounds on the second covariant derivatives of $\X$.

\begin{cor}     \label{cor:2nd}
Let $g\in\M^4(\S^n)$, and let $\X\colon(\S^n,g)\to\R^{n+1}$ be a $\C^4$
isometric embedding. Assume that $R$, the scalar curvature of $g$, is
positive. Then the following inequality holds:
\[
        \norm{\nabla^2 \X}^2 \leq
        \sup_{\S^n}\left( 2R - \frac{1}{R}\, \Delta R\right).
\]
\end{cor}
\begin{proof}
We have:
\begin{equation}        \label{eq:normal}
        \X_{;ij} = -\chi_{ij} \N,
\end{equation}
where $\N$ is the outer unit normal, and $\chi_{ij}$ the second
fundamental form of $\X$.  Thus, we obtain:
\[
        \norm{\nabla^2 \X}^2 = \X_{;ij} \cdot \X_;{}^{ij} = \chi^{ij}
        \chi_{ij}.
\]
Since $\chi^{ij}\chi_{ij} \leq H^2$, the corollary follows from
Theorem~\ref{thm:weyl}.
\end{proof}

We will use the Gauss and Codazzi Equations:
\begin{align}
        \label{eq:gauss}
        R_{ijkl} &= \chi_{ik} \chi_{jl} - \chi_{il} \chi_{jk} \\
        \label{eq:codazzi}
        0 &= \chi_{ij;k} - \chi_{ik;j},
\end{align}
where $R_{ijkl}$ is the Riemann curvature tensor of $g$. An immediate
consequence of~\eqref{eq:gauss} is that if $g$ has positive sectional
curvatures, then $\chi$ definite. In view of our choice of
normal~\eqref{eq:normal}, $\chi$ is positive definite.  Furthermore, it
follows from theorems of Hadamard, Chern and Lashof \cite{Chern}, and R.
Sacksteder \cite{Sa} that if $g$ has non-negative sectional curvature,
and $\X\colon(\S^n\,g)\to\R^{n+1}$ is an isometric immersion, then $\X$
is an embedding, and $\X(\S^n)$ is the boundary of a convex body in
$\R^{n+1}$.  In particular, if $g$ has non-negative sectional curvature,
then $\chi$ is positive semi-definite.

We begin the proof of Theorem~\ref{thm:weyl} with two lemmas.

\begin{lemma}   \label{lemma:sums}
Let $\x\in\R^n_+=\{(x_1,x_2,\dots,x_n)\in\R^n\st x_i>0\}$.  Then
\[
        3 \left(\sum x_i^2\right) \left(\sum x_i\right) \leq
        \left(\sum x_i\right)^3 + 2 \sum x_i^3.
\]
\end{lemma}

\label{page:sums}
\begin{proof}
This can be proved by induction.  However, we rather calculate directly:
\[
        \left(\sum x_i\right)^3 =
        \sum_{\norm{\alphab}=3} \binom{3}{\alphab} \, \x^{\alphab} =
        \sum_{\substack{\norm{\alphab}=3\\ \alphab^*=3}}
        \x^\alphab +
        3 \sum_{\substack{\norm{\alphab}=3\\ \alphab^*=2}}
        \x^\alphab +
        6 \sum_{\substack{\norm{\alphab}=3\\ \alphab^*=1}}
        \x^\alphab,
\]
Here $\alphab=(\alpha_1,\dots\alpha_n)$ stands for a multi-index, with
$\alpha_i$ nonnegative integers, $\norm{\alphab}=\sum \alpha_i$, and
$\alphab^*=\max \alpha_i$.  The multinomial coefficients are defined
for $\norm{\alphab}=k$ as
\[
        \binom{k}{\alphab}= \frac{k!}{\alphab!} =
        \frac{k!}{\alpha_1!\dots\alpha_n!},
\]
and $\x^\alphab=x_1^{\alpha_1} \dots x_n^{\alpha_n}$ for $\x\in\R^n$.
Furthermore, we have:
\[
        \left(\sum x_i^2\right) \left(\sum x_i\right) =
        \sum_{\substack{\norm{\alphab}=2,\>\alphab^*=2\\ \norm{\betab}=1}}
        \x^{\alphab+\betab} =
        \sum_{\substack{\norm{\alphab}=3\\ \alphab^*=3}}
        \x^\alphab +
        \sum_{\substack{\norm{\alphab}=3\\ \alphab^*=2}}
        \x^\alphab.
\]
Combining these, we conclude:
\[
        3 \left(\sum x_i^2\right) \left(\sum x_i\right) -
        \left(\sum x_i\right)^3 - 2 \sum x_i^3 =
        -6 \sum_{\substack{\norm{\alphab}=3\\ \alphab^*=1}}
        \x^\alphab \leq 0. \qed
\]
\renewcommand{\qed}{}
\end{proof}

Note that for $n=1,2$, the lemma becomes an identity.  As a corollary,
we obtain the following lemma.

\begin{lemma}   \label{lemma:chi}
Let $A$ be a symmetric, positive-definite $n\times n$ matrix.  Then we
have:
\[
        \left(\tr A\right)^3 - \tr A^3 \leq \frac32
        \left[ \left(\tr A\right)^2 - \tr A^2\right] \tr A.
\]
\end{lemma}
\begin{proof}
Let $x_1,x_2,\dots,x_n\geq0$ be the eigenvalues of $A$.  Then, by
Lemma~\ref{lemma:sums}, we have:
\begin{align*}
        2 \left[ \left(\tr A\right)^3 - \tr A^3 \right]
        &= 2 \left(\sum x_i \right)^3 - 2 \sum x_i^3 \\
        &\leq 3 \left[ \left(\sum x_i\right)^2 - \left(\sum x_i^2\right)
        \right] \sum x_i \\
        &= 3 \left[ \left(\tr A\right)^2 - \tr A^2 \right] \tr A.
        \qed
\end{align*}
\renewcommand{\qed}{}
\end{proof}

\begin{proof}[Proof of Theorem~\ref{thm:weyl}]
Suppose that $H$ achieves its maximum at $p\in\S^n$, then at $p$, we
have:
\[
        H_i=0, \qquad H_{;ij} \leq 0.
\]
Since $(\chi^{ij})\ge 0$ and $(H_{;ij})\le 0$, it follows with some
linear algebra that
\begin{equation}        \label{eq:deltah}
          H \Delta H \leq \chi^{ij} H_{;ij}
\end{equation}
Now, equation~\eqref{eq:gauss} implies that $R=H^2-\tr\chi^2$, where
$\tr\chi^2=\chi_{ij} \chi^{ij}$.  We also write
$\tr\chi^3=\chi^{ij}\chi_{jk}\chi^k{}_i$.  In view of~\eqref{eq:deltah},
\eqref{eq:codazzi}, \eqref{eq:gauss}, and Lemma~\ref{lemma:chi},
we now find that at $p$, the following holds:
\begin{equation}        \label{eq:long}
\begin{aligned}
        \Delta R &= 2 H \Delta H + 2 \norm{\nabla H}^2 -
        2\norm{\nabla\chi}^2 - 2\chi^{ij} \Delta \chi_{ij} \\
        &\leq 2 \chi^{ij} \left( \chi^k{}_{k;ji}
        - \chi_{ij;}{}^k{}_k \right) \\
        &= 2 \chi^{ij} \left( \chi^k{}_{j;ki} - \chi_j{}^k{}_{;ik}\right) \\
        &= 2 \chi^{ij} \left( R_{jlik} \chi^{lk} + R^k{}_{lik} \chi_j{}^l
        \right) \\
        &= 2 \chi^{ij} \left[ \left( \chi_{ji}\chi_{lk} - \chi_{jk}\chi_{li}
        \right) \chi^{lk} + \left( \chi^k{}_i\chi_{lk} -
        \chi^k{}_k\chi_{li}\right) \chi_j{}^l\right] \\
        &= 2 \left[ \left(\tr\chi^2\right)^2 - \tr\chi^3\tr\chi \right] \\
        &\leq 2 \left[ -2RH^2 + R^2 + \frac32 \left( \left(\tr\chi\right)^2
        -\tr\chi^2 \right) \left(\tr\chi\right)^2 \right] \\
        &= -RH^2 + 2R^2.
\end{aligned}
\end{equation}
It follows that
\[
        H^2 \leq H^2(p) \leq 2R(p) - \frac1{R(p)} \Delta R(p) \leq
        \sup_{\S^n} \left(2R - \frac1R \Delta R \right). \qed
\]
\renewcommand{\qed}{}
\end{proof}

\section{A Convergence Theorem}
\label{sec:converge}

In this section we show how to obtain $\C^{2,\alpha}$ a priori bounds
from Theorem~\ref{thm:weyl}.  Bounds of this type imply
Theorem~\ref{thm:convergence}.

Fix a finite covering $\{V_r\}$ of $\S^3$ by coordinate charts, and let
$\{U_r\}$ be a refinement such that $\overline U_r\subset V_r$.
Define for any $\C^{k,\alpha}$ tensor field $T$ of rank $l$ on $\S^3$
the norms:
\begin{align*}
        \Norm{T}_{k} &= \max_r \,
        \max_{\substack{0\leq\norm{\betab}\leq k \\ 1\leq
        i_1,\dots,i_l\leq n}}
        \sup_{x\in U_r}
        \norm{\D^\betab T_{i_1\dots i_l}(x)} \\[1ex]
        \left[ T \right]_{k,\alpha} &= \max_r
        \sup_{\substack{x,y\in U_r\\ \norm{\betab}=k}}
        \frac{\norm{\D^\betab T_{i_1\dots i_l}(x) - \D^\betab T_{i_1\dots
        i_l}(y)}}{\dist(x,y)^\alpha} \\[1ex]
        \Norm{T}_{k,\alpha} &=
        \max\left\{ \Norm{T}_{k}, \left[ T \right]_{k,\alpha} \right\}.
\end{align*}
where $T_{i_1\dots i_l}$ are the components of $T$ in the coordinate
chart on $U_r$. Endow the space of $\C^{k,\alpha}$ tensor fields of rank
$l$ over $\S^3$ with this norm.

\begin{thm}     \label{thm:holder}
Let $g\in\M^{k,\alpha}(\S^n)$ for some $k\geq4$ and $0<\alpha<1$, and
let $\X\colon(\S^n,g)\to\R^{n+1}$ be a $\C^{k,\alpha}$ isometric
embedding.  Suppose that the scalar curvature $R$ of $g$ is positive.
Then, there is $\X_0\in\R^{n+1}$, and a constant $C>0$ depending only on
$\Norm{g}_{k,\alpha}$, and $\min_{\S^n} R$ such that
\begin{equation}        \label{eq:holder}
        \Norm{\X-\X_0}_{k,\alpha}\leq C.
\end{equation}
\end{thm}
\begin{proof} It follows from the earlier mentioned theorems of
Hadamard, Chern and Lashof, and R. Sacksteder that $\X(\S^n)$ is the
boundary of an open convex set. $(\S^n,g)$ has bounded diameter which
implies that we have a bound on $\norm{\X-\X_0}$ provided that $\X_0$ is
chosen inside $\X(\S^n)$.  Next, from~\eqref{eq:iso}, we get that
$\norm{\nabla_i\X}^2=g_{ii}$ which clearly implies that
$\norm{\nabla\X}$ is bounded. Corollary~\ref{cor:2nd} implies bounds on
the second covariant derivatives of $\X$ depending only on $\Norm{g}_4$
and $\min_{\S^n}R$.  Since $\D_i\D_j \X = \nabla_i\nabla_j\X -
\Gamma^k_{ij} \nabla_i\X$, we get bounds on the second coordinate
derivatives of $\X$.  Next, we establish H\"older bounds on the second
derivatives of $\X$.  To proceed, we must show that $\X_0\in\R^{n+1}$
can be chosen so that there is a constant $C>0$, depending only on
$\Norm{g}_4$, such that
\begin{equation}        \label{eq:xdotn}
        (\X-\X_0)\cdot\N\geq\frac1C,
\end{equation}
where $\N$ is the outer unit normal to $\X$.

Let $K\subset\R^{n+1}$ be the closed convex set bounded by $\X(\S^n)$,
and choose $\X_0$ so that $B_r(\X_0)$ is a ball of largest radius
enclosed in $K$.  We claim that $r$ is bounded below by $1/C$, where $C$
is a constant which depends only on $\Norm{g}_4$.  Indeed, were this not
the case, then there would be a family
$g^\epsilon\in\M^{4,\alpha}(\S^n)$ with $\Norm{g^\epsilon}_4$ uniformly
bounded, and isometric embeddings
$\X^\epsilon\colon(\S^n,g^\epsilon)\to\R^{n+1}$ such that the largest
ball $B_{r^\epsilon}$ contained in the closed convex set $K^\epsilon$
bounded by $\X^\epsilon(\S^n)$ has $r^\epsilon\to0$, as $\epsilon\to0$.
This leads to a contradiction as follows.  By the Ascoli-Arzella
Theorem, there is a subsequence such that $g^{\epsilon_j}\to g$ in
$\C^{3,\beta}$ for any $0<\beta<1$.  Furthermore, by our argument above,
we have uniform bounds on $\Norm{\X^\epsilon}_2$, hence, perhaps along a
further subsequence, $\X^{\epsilon_j}\to\X$ in $\C^{1,\beta}$.  Note
that $\N^{\epsilon_j}$, the outer unit normals to $\X^{\epsilon_j}$,
converge in $\C^{\beta}$ to $\N$, the outer unit normal to $\X$. Indeed,
$\N^{\epsilon_j}$ is the exterior product of the $n$ vectors
$\X_{;1}^{\epsilon_j}, \dots, \X_{;n}^{\epsilon_j}$ divided by
$\det(g^{\epsilon_j})$.  Since, considering the normals as maps
$\N^{\epsilon_j}\colon\S^n\to\S^n$, we have $\deg(\N^{\epsilon_j})=1$,
we obtain $\deg(\N)=1$. Now, the limit of $K^{\epsilon_j}$ is a closed
convex set $K$ contained in a hyperplane. Otherwise, there would be
$n+1$ points in general position in $K$, which would imply that $K$
contained the simplex spanned by these $n+1$ points, which would
contradict $r^{\epsilon_j}\to 0$.  However, $\X(\S^n)$ is contained in
the same hyperplane, hence the image of $\N\colon\S^n\to\S^n$ consists
of at most 2 points, which contradicts $\deg(\N)=1$.
Inequality~\eqref{eq:xdotn} now follows easily, for if $\Pi_p$ is the
tangent plane  at $p\in\S^n$, then
\[
        \dist(p,\Pi_p) = (\X-\X_0)\cdot\N|_p.
\]
Since $K$ lies on one side of $\Pi_p$, and $B_{1/C}(\X_0)\subset K$,
we have $\dist(p,\Pi_p)\geq1/C$, and~\eqref{eq:xdotn} is established.
We now assume without loss of generality that $\X_0=\0$.

Define the function $\rho=\frac12\X\cdot\X$.  Then we have:
\begin{equation}        \label{eq:rhoij}
        \rho_{;ij} = g_{ij} - (\X\cdot\N)\, \chi_{ij},
\end{equation}
and consequently:
\begin{equation}        \label{eq:rho}
        \chi_{ij} = \frac{1}{\X\cdot\N}\, \bigl(g_{ij}-\rho_{;ij}\bigr).
\end{equation}
Let $\lambdab=(\lambda_1,\dots,\lambda_n)$ be the eigenvalues of
$g_{ij}-\rho_{;ij}$ (with respect to $g_{ij}$), and define for
$\x\in\R^n$:
\[
        f(\x) =
        \sum_{1\leq i<j \leq n} x_i x_j.
\]
The function $f$ is the second symmetric elementary function
$\sigma_{2,n}$; see Section~\ref{sec:gauss}.  It follows from
equation~\eqref{eq:rho} that
\begin{equation}        \label{eq:elliptic}
        f(\lambdab)^{1/2} = 2^{-1/2} (\X\cdot\N) R^{1/2}.
\end{equation}
Let $\Gammab_2$ denote the connected component of the set
$\{\x\in\R^n\st f(\x)>0\}$ containing the positive cone.  It follows
from~\cite{caffarelli85} that $\Gammab_2$ is a cone with the property
that:
\begin{align*}
        \frac{\D}{\D x_i}\left( f(\x)^{1/2}
        \right) &> 0, \quad \forall\x\in\Gammab_2, \>
                \forall 1\leq i\leq n; \\
        \left( \frac{\D^2}{\D x_i\D x_j}\left(
        f(\x)^{1/2} \right)\right) &\leq 0, \quad
        \forall\x\in\Gammab_2,
\end{align*}
see also~\cite{garding59}. Now, in view of our $\C^0$ bound,
\eqref{eq:xdotn}, and our hypothesis on $R$, we have an estimate:
\[
        \frac1C\leq (\X\cdot\N) R^{1/2} \leq C,
\]
hence the eigenvalues $\lambdab=(\lambda_1,\dots,\lambda_n)$  of
$g_{ij}-\rho_{;ij}$ with respect to $g_{ij}$ remain within a fixed
compact set of $\Gammab_2$. We conclude that
Equation~\eqref{eq:elliptic}, as an equation in the Hessian $\rho_{;ij}$
of $\rho$, is uniformly elliptic and convex in $\rho_{;ij}$. As we
already have estimates on $\Norm{\rho}_2$, it now follows
from~\cite[Theorem 6.6 and Theorem 8.1]{caffarelli95} that there is a
constant $C$ depending only on $\Norm{g}_{4,\alpha}$ and $\min_{\S^2} R$
such that:
\begin{equation}        \label{eq:rho2a}
        \Norm{\rho}_{2,\alpha} \leq C;
\end{equation}
see also~\cite{evans82} and~\cite{krylov82}. Once we have
established~\eqref{eq:rho2a}, we can apply Schauder estimates
to~\eqref{eq:elliptic} to get
\begin{equation}        \label{eq:rho4a}
        \Norm{\rho}_{k,\alpha} \leq C.
\end{equation}
In view of~\eqref{eq:normal} and~\eqref{eq:rho}, we have:
\[
        \X_{;ij} = -\frac{1}{\X\cdot\N}\, \bigl(g_{ij}-\rho_{;ij}\bigr)
        \N.
\]
Thus, \eqref{eq:rho4a} implies~\eqref{eq:holder}.
\end{proof}

\begin{proof}[Proof of Theorem~\ref{thm:convergence}]
Follows directly from Theorem~\ref{thm:holder}.
\end{proof}

\section{Weyl-Type Estimates for Non-Negative Scalar Curvature}
\label{sec:weyltype}

In this section, we prove Theorem~\ref{thm:guanli94}. The proof of
Theorem~\ref{thm:guanli94} is similar to that of the Proposition
in~\cite{guanli94}, using~\eqref{eq:long} instead of the Weyl estimate.
We first point out a consequence of Theorem~\ref{thm:guanli94} which is
similar to Corollary~\ref{cor:2nd}.

\begin{cor}     \label{cor:3rd}
Let $g\in\M^4(\S^n)$, and
let $\X\colon(\S^n,g)\to\R^{n+1}$ be a $\C^4$ isometric embedding.
Let $R$ be the scalar curvature of $g$.
Then the following inequality holds:
\[
        \norm{\nabla^2 \X}^2 \leq C d^2
        \sup_{\S^n} \left(2R^2-\Delta R
        + \frac{(n-1)^2}{64 d^2} R\right),
\]
where $C= 4(n-1)^{-2}e^{(n-1)/4}$, and $d$ is the diameter of
$(\S^n,g)$.
\end{cor}

The proof of Corollary~\ref{cor:3rd} is the same as that of
Corollary~\ref{cor:2nd}, using  Theorem~\ref{thm:guanli94}
instead of Theorem~\ref{thm:weyl}.

\begin{proof}[Proof of Theorem~\ref{thm:guanli94}]
Define $f=e^{\alpha\rho} H$, where $\alpha>0$ is a constant to be
determined later, and as before $\rho=\frac12\X\cdot\X$.  Suppose that
$f$ achieves it maximum at $p\in\S^n$, then at $p$, we have:
\begin{align}
        \label{eq:fi}
        f_i &= \alpha \rho_i H + H_i = 0, \\
        \label{eq:fij}
        f_{;ij} &= e^{\alpha\rho} (\alpha H \rho_{;ij} - \alpha^2 H
        \rho_i\rho_j + H_{;ij}) \leq 0.
\end{align}
Furthermore, note that
\[
        \X = g^{ij} \rho_i \X_j + (\X\cdot\N)\, \N.
\]
This is verified by taking inner product with the $n+1$ independent
vectors $\X_{;1},\dots,\X_{;n},\N\in\R^{n+1}$.
On the other hand, taking inner product with $\X$ we obtain:
\begin{equation}        \label{eq:gradrho}
        2\rho = \norm{\nabla\rho}^2 + (\X\cdot\N)^2.
\end{equation}
In particular, since $2\rho\leq d^2$, we have the following inequalities:
\begin{align}
        \label{eq:estgradrho}
        \norm{\nabla\rho}^2 &\leq d^2, \\
        \label{eq:estXdotN}
        \X\cdot\N &\leq d.
\end{align}
As in~\eqref{eq:long}, we now calculate at $p$, taking
~\eqref{eq:fi} and ~\eqref{eq:estgradrho} into account:
\begin{equation}        \label{eq:longer}
\begin{aligned}
        \Delta R &= 2 H \Delta H + 2 \norm{\nabla H}^2 -
        2\norm{\nabla\chi}^2 - 2\chi^{ij} \Delta \chi_{ij} \\
        &\leq 2H\Delta H - 2 \chi^{ij} H_{;ij}
        + 2 \alpha^2 \norm{\nabla\rho}^2 H^2
        + 2 \chi^{ij} \left( R_{jlik} \chi^{lk} + R^k{}_{lik}
        \chi_j{}^l \right) \\
        &\leq 2(Hg^{ij} - \chi^{ij}) H_{;ij} + 2 \alpha^2 d^2 H^2
        - RH^2 + 2R^2.
\end{aligned}
\end{equation}
Since $(\chi^{ij})\geq0$ and $(f_{;ij})\le 0$ at $p$,
we have,
$$
\chi^{ij}f_{;ij}\ge Hg^{ij}f_{;ij},
$$
namely,
$$
(H g^{ij} - \chi^{ij})f_{;ij}\le 0.
$$
Thus,
 in view of~\eqref{eq:fij}, \eqref{eq:rhoij},
and~\eqref{eq:estXdotN}:
\begin{equation}        \label{eq:DeltaH}
\begin{aligned}
        (H g^{ij} - \chi^{ij}) H_{;ij}
        &\leq (Hg^{ij}-\chi^{ij}) (-\alpha H\rho_{ij} + \alpha^2
        H\rho_i\rho_j) \\
        &\leq -\alpha H(Hg^{ij}-\chi^{ij})(g_{ij} - (\X\cdot\N)\,\chi_{ij})
        + \alpha^2 H^2 \norm{\nabla\rho}^2 \\
        &\leq -(n-1)\alpha H^2 + \alpha RH\, (\X\cdot\N) +
        \alpha^2 d^2 H^2 \\
        &\leq -(n-1)\alpha H^2 + \alpha^2 d^2 H^2 + \alpha d RH.
\end{aligned}
\end{equation}
Thus, combining~\eqref{eq:longer} and~\eqref{eq:DeltaH}, we obtain that
the following inequality holds at $p$:
\[
        \Delta R \leq -2(n-1)\alpha H^2 + 4\alpha^2 d^2 H^2 +
        RH(\alpha d - H) + 2R^2.
\]
Taking $\alpha=(n-1)/(4d^2)$, and using $H(\alpha d-H) \leq \alpha^2
d^2/4$, we get that the following inequality holds at $p$:
\[
        \Delta R \leq - \frac{(n-1)^2}{4d^2} H^2 + 2R^2 +
        \frac{(n-1)^2}{64 d^2} R.
\]
Thus, we conclude
\begin{equation}        \label{eq:H(p)}
        H^2(p) \leq \frac{4 d^2}{(n-1)^2}\, \sup_{\S^n} \left(
        2R^2 - \Delta R + \frac{(n-1)^2}{64 d^2} R \right).
\end{equation}
If $q$ is any point on $\S^n$, then we conclude:
\begin{multline}        \label{eq:H(q)}
        H^2(q) =e^{-2\alpha\rho(q)} f^2(q)
        \leq e^{-2\alpha\rho(q)} f^2(p) \\
        = e^{2\alpha(\rho(p)-\rho(q))} H^2(p)
        \leq e^{\alpha d^2} H^2(p)
        \leq e^{(n-1)/4} H^2(p).
\end{multline}
Theorem~\ref{thm:guanli94} now
follows from~\eqref{eq:H(p)} and~\eqref{eq:H(q)}.
\end{proof}

\section{Solving the Once-Contracted Gauss Equation}
\label{sec:gauss}

In this section, we show that when the sectional curvature of $g$ is
strictly positive, one can bound $H$ in terms of the supremum of the
Ricci tensor and the minimum of the sectional curvature, \ie in terms of
only two derivatives of $g$.  This is done by showing that the
once-contracted Gauss Equation:
\begin{equation}        \label{eq:tracegauss}
        R_{ij} = \tr\chi\, \chi_{ij} - \chi_i{}^k \chi_{kj},
\end{equation}
can be solved for the second fundamental form $\chi$;
cf.~Theorem~\ref{thm:solve}.  It is well known that the full Gauss
equation has at most one solution when $n\geq3$ and the sectional
curvature is positive, see for example~\cite[Section 60]{eisenhart26}
and~\cite{allendoerfer39}. Here we show that the subset consisting
of~\eqref{eq:tracegauss} can always be solved under this assumption. As
a consequence, we obtain an estimate on the second derivatives of $\X$
which depends only on the second derivatives of the metric $g$, provided
the sectional curvatures of $g$ are positive;
cf.~Theorem~\ref{thm:c2bound}.  When $n=3$, this also gives, in
conjunction with the Codazzi equation~\eqref{eq:codazzi}, an explicit
local necessary and sufficient intrinsic condition for a metric $g$ of
positive sectional curvature on $\S^3$ to be locally embeddable in
$\R^{4}$, and consequently also globally embeddable;
cf.~Theorem~\ref{thm:localembed} and the remark following it.

We will use the elementary symmetric functions $\sigma_{k,n}(\x)$
defined for integers $k\geq0$ and $n\geq1$, and for
$\x=(x_1,\dots,x_n)\in\R^n$, by:
\begin{gather*}
        \sigma_{0,n}(\x) = 1, \quad
        \sigma_{k,n}(\x) = 0, \>\text{if $k>n$}, \\
        \sigma_{k,n}(\x) = \sum_{i_1<\dots<i_k} x_{i_1} \dots x_{i_k},
        \>\text{for $1\leq k\leq n$}.
\end{gather*}
Let $\x_i=(x_1,\dots,\hat x_i,\dots,x_n)\in\R^{n-1}$ be
the vector obtained from $\x$ by deleting its $i$-th coordinate.
We have for $k\geq1$, and $1\leq i\leq n$:
\begin{equation}        \label{eq:recursive}
        \sigma_{k,n}(\x) = x_i \sigma_{k-1,n-1}(\x_i)
        + \sigma_{k,n-1}(\x_i),
\end{equation}
from which it follows by induction on $n$ that:
\begin{equation} \label{eq:sigma}
        \sum_{i=1}^n \sigma_{k,n-1}(\x_i) = (n-k) \sigma_{k,n}(\x),
\end{equation}
for $k\geq0$.

\begin{lemma}   \label{lemma:det(g)}
Let $n\geq3$, $\x=(x_1,\dots,,x_n)\in\R^n$, $x=\sigma_{1,n}(\x)$, and
define the $n\times n$ matrix $G_n(\x)$ by:
\[
        G_n(\x) = \begin{pmatrix}
        x-x_1 & x_1 & \dots & x_1 \\
        x_2 & x-x_2 & \dots & x_2 \\[.5ex]
        \vdots \\[.5ex]
        x_n & x_n & \dots & x-x_n
        \end{pmatrix}.
\]
Then:
\[
        \det G_n(\x)= \sum_{\norm{\gammab}=n} a_{\gammab,n} \x^\gammab,
\]
where
\begin{equation}        \label{eq:det(g)}
        a_{\gammab,n} = \sum_{k=3}^n \frac{(-2)^{k-1}(k-2)
        (n-k)!}{\gammab!} \sigma_{k,n}(\gammab).
\end{equation}
\end{lemma}
\begin{proof}
Define for $\x\in\R^n$, and $s\in\R$:
\begin{align*}
        F_n(s,\x) &= \det \begin{pmatrix}
        s-x_1 & x_1 & \dots & x_1 \\
        x_2 & s-x_2 & \dots & x_2 \\[.5ex]
        \vdots \\[.5ex]
        x_n & x_n & \dots & s-x_n
        \end{pmatrix}, \\[1ex]
        f_n(s,\x) &= \det F_n(s,\x).
\end{align*}
We will show that:
\begin{equation}        \label{eq:fn}
        f_n(s,\x) = s^n - \sigma_{1,n}(\x) s^{n-1} + \sum_{k=3}^n (-2)^{k-1}
        (k-2) \sigma_{k,n}(\x) s^{n-k}.
\end{equation}
This is shown by induction on $n\geq3$.  Equation~\eqref{eq:fn} is
easily verified for $n=3$.  Assume~\eqref{eq:fn} holds for $n-1\geq3$,
and let $g_n(s,\x)$ denote the right hand side of~\eqref{eq:fn}. Now,
using
\[
        \frac{d}{ds} \log\det F_n(s,\x) =
        \tr\left(F_n(s,\x)^{-1} \frac{d}{ds} F_n(s,\x) \right),
\]
the induction hypothesis, and~\eqref{eq:sigma}, we obtain:
\[
        \frac{d}{ds} f_n(s,\x) = \sum_{i=1}^n f_{n-1}(s,\x_i) =
        \sum_{i=1}^n g_{n-1}(s,\x_i) = \frac{d}{ds} g_n(s,\x).
\]
It is easily checked that $f_n(0,\x)=(-2)^{n-1}(n-2)\sigma_{n,n}(\x)
=g_n(0,\x)$, hence $f_n(s,\x)=g_n(s,\x)$, which
proves~\eqref{eq:fn}.  Substituting $s=\sigma_{1,n}(\x)$
in~\eqref{eq:fn}, we now obtain:
\[
        \det G_n(\x) = \sum_{k=3}^n (-2)^{k-1}
        (k-2) \sigma_{k,n}(\x) \sigma_{1,n}(\x)^{n-k}.
\]
We have
\[
        \sigma_{1,n}(\x)^{n-k} = \sum_{\norm{\alphab}=n-k}
        \binom{n-k}{\alphab} \x^\alphab, \qquad
        \sigma_{k,n}(\x) = \sum_{\substack{\norm{\betab}=k\\
        \betab^*=1}} \x^\beta,
\]
see the proof of Lemma~\ref{lemma:sums} on page~\pageref{page:sums}.
Thus, substituting $\gammab=\alphab+\betab$, we find
\[
        \det G_n(\x) = \sum_{k=3}^n (-2)^{k-1} (k-2)
        \sum_{\norm{\gammab}=n} \,
        \sum_{\substack{\norm{\betab}=k \\
        \betab\leq\gammab,\>\betab^*=1}}
        \binom{n-k}{\gammab-\betab} \x^\gammab
        = \sum_{\norm{\gammab}=n} a_{\gammab,n} \x^\gammab,
\]
where
\[
        a_{\gammab,n} = \sum_{k=3}^n (-2)^{k-1} (k-2)
        \sum_{\substack{\norm{\betab}=k \\
        \betab\leq\gammab,\> \betab^*=1}} \binom{n-
        k}{\gammab-\betab}.
\]
To obtain~\eqref{eq:det(g)}, note
that the last sum in this equation can be rewritten as:
\[
        \sum_{\substack{\norm{\betab}=k \\
        \betab\leq\gammab,\> \betab^*=1}} \binom{n-
        k}{\gammab-\betab} =
        \frac{(n-k)!}{\gammab!}\sigma_{k,n}(\gammab).
\]
\end{proof}

\begin{lemma}   \label{lemma:det}
Let $n\geq3$, and $\x\in\R^n_+$, then $\det G_n(\x)>0$.
\end{lemma}
\begin{proof}
Let
\begin{equation}        \label{eq:defb}
        b_{\gammab,k,n} = \frac{2^{k-1}(k-2)(n-k)!}{\gammab!}
        \sigma_{k,n}(\gammab).
\end{equation}
Then, we have:
\begin{equation}        \label{eq:a}
        a_{\gammab,n} = \sum_{k=3}^n (-1)^{k-1} b_{\gammab,k,n}.
\end{equation}
We now claim that for all $\norm{\gammab}=n$, and all $k\geq 3$, there
holds:
\begin{equation}        \label{eq:b}
        2(k-1)\sigma_{k+1,n}(\gammab) \leq (n-k)(k-2)\sigma_{k,n}(\gammab),
\end{equation}
It follows from~\eqref{eq:b} and~\eqref{eq:defb}
that $b_{\gammab,k+1,n}\leq b_{\gammab,k,n} $,
for all $\norm{\gammab}=n$, and $3\leq k\leq n-1$, and hence, in view
of~\eqref{eq:a}, we obtain $a_{\gammab,n}\geq0$, for all $\norm{\gammab}=n$.
For $n=3$, we have $\det G(x_1,x_2,x_3)= 4x_1 x_2 x_3$, that is
$a_{\gammab,3}=0$ unless $\gammab=(1,1,1)$, in which case $a_{\gammab,3}=4$.
Since for $n>3$, we have
\[
        \det G(x_1,x_2,x_3,0,\dots,0)=(x_1+x_2+x_3)^{n-3}\det
        G(x_1,x_2,x_3),
\]
it is clearly impossible that $a_{\gammab,n}=0$ for all
$\norm{\gammab}=n$, and the lemma follows.  It remains to
prove~\eqref{eq:b}.  In fact, we will prove the more general inequality:
\begin{equation}        \label{eq:ineq}
        (k+1)\sigma_{k+1,n}(\gammab) \leq
        \bigl(\sigma_{1,n}(\gammab)-k\bigr)\sigma_{k,n}(\gammab),
\end{equation}
for $k\geq0$, and any multi-index $\gammab$.
This clearly implies~\eqref{eq:b} when $\norm{\gammab}=n$, and
$k\geq3$.    We will use the identity:
\[
        \sum_{i=1}^n \gamma_i\sigma_{k,n-1}(\gammab_i) =
        (k+1) \sigma_{k+1,n}(\gammab),
\]
which holds for $k\geq0$, and
which follows easily from~\eqref{eq:recursive} by induction on $n$.
Using this identity, we find that:
\begin{align*}
        \sigma_{1,n}(\gammab) \sigma_{k.n}(\gammab) &=
        \sum_{i=1}^n \gamma_i \bigl(\gamma_i \sigma_{k-1}(\gammab_i) +
        \sigma_{k,n-1}(\gammab_i)\bigr) \\
        &\geq \sum_{i=1}^n \bigl( \gamma_i \sigma_{k,n}(\gammab_i)
        + \gammab_i \sigma_{k,n-1}(\gammab_i)\bigr) \\[1ex]
        & = (k+1) \sigma_{k+1,n}(\gammab) + k \sigma_{k,n}(\gammab),
\end{align*}
which proves~\eqref{eq:ineq}.
\end{proof}

\begin{lemma}   \label{lemma:1-1}
Let $n \geq 3$, and let $\Sn_n$ be the set of positive-definite
symmetric $n\times n$ matrices over $\R$.  The map $\Phi\colon \Sn_n\to
\Sn_n$, defined by
\[
        \Phi\colon A \mapsto (\tr A) A - A^2
\]
is one-to-one.
\end{lemma}
\begin{proof}
Suppose that $\Phi(A)=\Phi(B)$, and let $E=A+B$, $F=A-B$.  Then,
$E$ and $F$ are symmetric, and $E$ is positive-definite.  We have:
\begin{equation}        \label{eq:uniqueness}
        \tr(E) F + \tr(F) E - (EF+FE) = 0.
\end{equation}
This implies that $E$ and $F$ can be simultaneously diagonalized.
Indeed, let $\lambdab=(\lambda_1,\dots,\lambda_n)$ be the eigenvalues of
$E$, and let $\mub=(\mu_1,\dots,\mu_n)$ be the eigenvalues of $F$, then
there is an orthogonal matrix $Q$ such that
$Q^TEQ=D=\diag(\lambdab)$.  Let $C=Q^TFQ$, then we have
\[
        \tr(D)C+\tr(C)D-(DC+CD)=0.
\]
Writing $C\e_j=\sum_{i=1}^n c_{ij}\e_i$, where $\e_1,\dots,\e_n$ is the
standard basis of $\R^n$, this implies:
\[
        \sum_{i=1}^n
        \bigl(\sigma_{1,n}(\lambdab)-(\lambda_i+\lambda_j)\bigr)
        c_{ij}\e_i = -\lambda_j\sigma_{1,n}(\mub)\e_j.
\]
Since $\sigma_{1,n}(\lambdab)-(\lambda_i+\lambda_j)>0$ when $i\ne j$,
this shows that $C$ is diagonal. It now follows
from~\eqref{eq:uniqueness} that
\[
        \mu_i\sum_{j\ne i} \lambda_j + \lambda_i\sum_{j\ne i} \mu_j = 0.
\]
In view of Lemma~\ref{lemma:det}, we obtain that $\mu_i=0$ for
$i=1,\dots,n$.  Thus, we conclude that $F=0$, and $A=B$.
\end{proof}

Denote by $\Tn_n\subset\Sn_n$ the set of those matrices $B\in \Sn_n$
with the following property: \emph{if $\mu_1,\dots,\mu_n$ are the
eigenvalues of $B$, then $\mu_i<\sum_{j\ne i}\mu_j$ for each $1\leq
i\leq n$.}  For $A\in\Sn_n$, let $\norm{A}^2=\sum_{i=1}^n
\lambda_i^2$, where $\lambda_i$ are the eigenvalues of $A$,
and for $B\in\Tn_n$, let
\[
        \varepsilon(B) = \min_{1\leq i\leq n} \bigl(\sigma_{1,n}(\mub)-
        2\mu_i)>0,
\]
where $\mub=(\mu_1,\dots,\mu_n)$ are the eigenvalues of $B$.

\begin{lemma}   \label{lemma:onto}
Let $n\geq3$, then
$\Phi$ maps $\Sn_n$ onto $\Tn_n$.  Furthermore if
$A\in\Sn_n$, then
\begin{equation}        \label{eq:normA}
        \norm{A} \leq \frac n2 \norm{\Phi(A)}
        \varepsilon\bigl(\Phi(A)\bigr)^{-1/2}.
\end{equation}
\end{lemma}
\begin{proof}
We first prove that $\Phi(\Sn_n)\subset\Tn_n$.  Let $A\in\Sn_n$, and let
$\lambdab=(\lambda_1\dots,\lambda_n)$ be its eigenvalues. Then, in view
of~\eqref{eq:recursive}, the eigenvalues of $\Phi(A)$ are:
\[
        \mu_i = \sum_{j\ne i} \lambda_j\lambda_i =
        \lambda_i \sigma_{1,n-1}(\lambdab_i)
        = \sigma_{2,n}(\lambdab) - \sigma_{2,n-1}(\lambdab_i).
\]
Thus, using~\eqref{eq:sigma}, we find
that $\sigma_{1,n}(\mub)=2\sigma_{2,n}(\lambdab)$, and
hence:
\begin{equation}        \label{e(B)}
        \sigma_{1,n}(\mub) - 2\mu_i = 2\sigma_{2,n-1}(\lambdab_i) > 0.
\end{equation}
We conclude that $\Phi(A)\in\Tn_n$.  The converse will be proved by
continuity, \ie we will show that $\Phi(\Sn_n)$ is open and closed in
$\Tn_n$.  Since $\Tn_n$ is clearly connected, this implies that
$\Phi(\Sn_n)=\Tn_n$.  To show that $\Phi(\Sn_n)$ is open suppose that
$A\in\Sn_n$, and that
\begin{equation}        \label{eq:linearized}
        \Phi'(A)C=\tr(A)C + \tr(C)A - (AC+CA) = 0,
\end{equation}
for some symmetric $n\times n$ matrix $C$.  Then, as in the proof of
Lemma~\ref{lemma:1-1}, we get that $C=0$, and hence $\Phi'(A)$ is
non-singular.  We conclude, by the Inverse Function Theorem, that there
is a neighborhood of $\Phi(A)$ contained in $\Phi(\Sn_n)$, and thus,
$\Phi(\Sn_n)$ is open.  We will now prove~\eqref{eq:normA}.  Let
$A\in\Sn_n$, and let $\lambdab=(\lambda_1,\dots,\lambda_n)$ be its
eigenvalues. Let $\lambda_i=\max_j\lambda_j$, so that $\norm{A}^2 \leq
n\lambda_i^2$. Then, in view of~\eqref{eq:recursive} and~\eqref{e(B)},
we have:
\[
        \sigma_{2,n}(\lambdab)^2 \geq
        \lambda_i^2 \sigma_{1,n-1}(\lambdab_i)^2 \geq
        \frac2n \norm{A}^2 \sigma_{2,n-1}(\lambdab_i) \geq
        \frac1n \norm{A}^2 \varepsilon\bigl(\Phi(A)\bigr).
\]
If $B=\Phi(A)$, we can now estimate:
\[
        \norm{B}^2 \geq \frac1n (\tr B)^2 =
        \frac1n \bigl((\tr A)^2 - \tr(A^2)\bigr)^2
        = \frac4n \sigma_{2,n}(\lambdab)^2 \geq
        \frac4{n ^2} \norm{A}^2 \varepsilon\bigl(\Phi(A)\bigr),
\]
which proves~\eqref{eq:normA}.  This shows that $\Phi(\Sn_n)$ is closed
in $\Tn_n$.  Indeed, if $A_j\in\Sn_n$ is such that $\Phi(A_j)\to
B\in\Tn_n$, then~\eqref{eq:normA} shows that the eigenvalues
$\lambdab(A_j)$ of $A_j$ are uniformly bounded above.  Therefore, by
passing to a subsequence if necessary, we see that $A_j\to A$ for some
symmetric $A\geq0$.  By continuity $\Phi(A)=B$, hence $A>0$, \ie
$A\in\Sn_n$, and $B\in\Phi(\Sn_n)$.
\end{proof}

We note here that when $n=3$, the inverse of $\Phi$ can be written
explicitly.  Indeed, let $B\in\Tn_3$, let $\mub=(\mu_1,\mu_2,\mu_3)$ be
its eigenvalues, and let $\v_1,\v_2,\v_3$ be its eigenvectors.  Define:
\begin{equation}        \label{eq:case3}
        \lambda_i = \frac{\sqrt{\prod_{j=1}^3 \bigl(\sigma_{1,3}(\mub)-
        2\mu_j\bigr)}}{\sqrt{2}\bigl(\sigma_{1,3}(\mub)-2\mu_i\bigr)},
        \> i=1,2,3;
        \quad A = \sum_{i=1}^3 \lambda_i \v_i\otimes\v_i.
\end{equation}
It is easy to check that $\lambda_i\sum_{j\ne i}\lambda_j=\mu_i$, hence
$\Phi(A)=B$.

We note, furthermore, that since $\Phi\colon\Sn_n\to\Tn_n$ is
$\C^\infty$ (in fact analytic), it follows that
$\Phi^{-1}\colon\Tn_n\to\Sn_n$ is also of class $\C^\infty$.

\begin{thm}     \label{thm:solve}
Let $n\geq3$, and let $g\in \M^2_+(\S^n)$.  Then there exists a unique
symmetric positive twice-covariant tensor $\chi$ on $\S^n$ which
satisfies the once-contracted Gauss Equations~\eqref{eq:tracegauss}.
Furthermore, if $g\in\C^{k,\alpha}$ for some $k\geq2$ and
$0\leq\alpha<1$, then $\chi\in\C^{k-2,\alpha}$.
\end{thm}
\begin{proof}
Let $R_{ij}$ be the Ricci tensor of $g$ in an orthonormal basis, and let
$p\in\S^n$.  Then, we have $R_{ij}(p)\in\Tn_n$.  Indeed, let $\mu_i$ be
the eigenvalues of $R_{ij}$, let $\kappa_{ij}$ be the sectional
curvature of $g$ in the plane spanned by the $i$-th and $j$-th
eigenvector of $R_{ij}$, and write $\kappa_{ii}=0$.  Then,
$\mu_i=\sum_{j=1}^n \kappa_{ij}$ for each $i$, hence:
\begin{equation}        \label{eq:curvatures}
        \sum_{j\ne i} \mu_{j} =
        \sum_{j\ne i} \sum_{l=1}^n \kappa_{jl}
        = \sum_{j\ne i} \sum_{l\ne i} \kappa_{jl} + \mu_i > \mu_i.
\end{equation}
Thus, the first assertion of the theorem follows from
Lemma~\ref{lemma:onto} and Lemma~\ref{lemma:1-1}.  The second assertion
follows from the remark just preceding this theorem.
\end{proof}

This theorem, in conjunction with~\eqref{eq:normA} and~\eqref{eq:curvatures},
immediately implies the following theorem:

\begin{thm}     \label{thm:c2bound}
Let $n\geq3$, let $g\in\M^2_+(\S^n)$, and let
$\X\colon(\S^n,g)\to\R^{n+1}$ be a $\C^2$ isometric embedding.  Let
$\chi$ be the second fundamental form of $\X$, let $\kappa$ be the
minimum of all the sectional curvatures of $g$ on $\S^n$, and let
$\Lambda$ be the maximum of the norm
$\norm{\operatorname{Ric}}=(R_{ij}R^{ij})^{1/2}$ of the Ricci tensor of
$g$ over $\S^n$.  Then, we have:
\begin{equation}        \label{eq:c2bound}
        \norm{\chi}\leq C_n\Lambda \kappa^{-1/2},
\end{equation}
where $C_n=n/\bigl(2\sqrt{(n-1)(n-2)}\bigr)$.
\end{thm}

This last theorem should be compared with Theorems~\ref{thm:weyl}
and~\ref{thm:guanli94} when $n\geq3$.  Here we require positive
sectional curvature. However, we obtain an estimate which depends only
on two derivatives of $g$. In Theorem~\ref{thm:weyl}, we only required
non-negative sectional curvature and positive scalar curvature, but our
estimate relied on four derivatives of $g$.

\begin{defn}
Let $(M,g)$ be a Riemannian manifold of class $\C^{k,\alpha}$ for some
$k\geq2$ and $0\leq\alpha<1$. We say that $(M,g)$ is
\emph{$\C^{k,\alpha}$ locally isometrically embeddable} in $\R^N$ if for
each $p\in M$ there is a neighborhood $U\subset M$ of $p$ and a
$\C^{k,\alpha}$ isometric embedding $\X\colon (U,g)\to\R^N$.  If $(M,g)$
is locally isometrically embeddable in $\R^N$, we say that $(M,g)$ is
\emph{locally rigid} in $\R^N$ if whenever $\X,\X'\colon (U,g)\to\R^N$
are local $\C^2$ isometric embeddings of some open set $U\subset M$,
then $\X'=\Psi\compose\X$, where $\Psi\colon\R^N\to\R^N$ is a rigid
motion possibly composed with a reflection.
\end{defn}

We now have:

\begin{thm}     \label{thm:localembed}
Let $g\in\M^{k,\alpha}_+(\S^3)$ for some $k\geq3$ and $0\leq\alpha<1$.
Let $\chi$ be the solution of the
once-contracted Gauss Equation~\eqref{eq:tracegauss} given by
Theorem~\ref{thm:solve}.  Then $(\S^3,g)$ is $\C^{k,\alpha}$
locally isometrically embeddable in
$\R^{4}$ if and only if $\chi$ satisfies the Codazzi
Equations:
\begin{equation}        \label{eq:intrinsic}
        \chi_{ij;k} - \chi_{ik;j}=0
\end{equation}
In this case, $(\S^3,g)$ is also locally rigid.
\end{thm}
\begin{proof}
Let $\X\colon(U,g)\to\R^4$ be a $\C^{k,\alpha}$ local isometric
embedding.  Then, by Theorem~\ref{thm:solve}, $\chi$ is the second
fundamental form of $\X$, and $\chi\in\C^{k-2,\alpha}$. Hence $\chi$
satisfies the Codazzi Equations.  Conversely, the system
\begin{equation}        \label{eq:local}
        \left\{
        \begin{aligned}
        \X_{;ij} &= -\chi_{ij} \N \\
        \N_i &= \chi_i{}^j \X_{;j}
        \end{aligned}
        \right.
\end{equation}
is an overdetermined system for the $4$ vector fields
$\X_{;1},\X_{;2},\X_{;3},\N$ along $\S^3$, whose integrability
conditions are the Gauss and Codazzi
Equations~\eqref{eq:gauss}--\eqref{eq:codazzi}.  The Codazzi Equations
are satisfied by hypothesis, and the once-contracted Gauss
Equations~\eqref{eq:tracegauss} imply the full Gauss
Equations~\eqref{eq:gauss} when $n=3$. Thus, if $p\in\S^3$, we can
integrate~\eqref{eq:local} in a neighborhood $U$ of $p$ in $\S^3$ with
$\X_{;i}\cdot\X_{;j}|_p=g_{ij}(p)$.  Furthermore, since
$\X_{;ij}=\X_{;ji}$ we can integrate again to obtain $\X$.  It then
follows that $\X_{;i}\cdot\X_{;j}=g_{ij}$ throughout $U$.  Since
$\chi\in\C^{k-2,\alpha}$, we obtain from~\eqref{eq:local} that
$\X\in\C^{k,\alpha}$.  The rigidity statement follows from the
uniqueness of $\chi$.
\end{proof}

\begin{rem}
As mentioned in the introduction, local embeddability implies global
embeddability under the hypothesis of positive sectional curvature.
 Therefore as a corollary of the previous theorem we
can replace \emph{local} by \emph{global} in the previous theorem.
\end{rem}

\end{document}